\documentclass[aos,MSNbibl,seceqn,nameyear,dvips]{arximspdf}
\usepackage{mathrsfs}

\doi{10.1214/13-AOS1098} 
\volume{41}
\issue{3}
\pubyear{2013}
\firstpage{1142}
\lastpage{1165}

\makeatletter

\newtheorem{theorem}{Theorem}[section]
\newtheorem{lemma}{Lemma}[section]
\newtheorem{corollary}{Corollary}[section]

\newcommand{\cal}{\mathcal}

\def\bb{{\mathbf b}}
\def\bG{{\mathbf G}}
\def\dell{\dot{\ell}}
\def\ddell{\ddot{\ell}}
\def\Nbar{\bar{N}}
\def\cO{{\cal O}}
\def\rP{{\mathrm P}}
\def\bu{{\mathbf u}}
\def\tx{\tilde{x}}
\def\bX{{\mathbf X}}
\def\by{{\mathbf y}}
\def\bZ{{\mathbf Z}}

\def\eps{\varepsilon}
\def\argmin{\mathop{\arg\min}}
\def\real{\mathbb R}

\def\sgn{\operatorname{sgn}}

\def\hbeta{\hat{\beta}}

\def\bu{{\mathbf u}}
\def\bX{{\mathbf X}}
\def\btheta{\bolds\theta}
\def\bbeta{\bolds\beta}

\def\bhbeta{\hat{\bbeta}}

\def\hbbeta{\hat{\bolds\beta}}
\def\lam{\lambda}
\def\ttheta{\tilde{\theta}}

\def\tbtheta{\tilde{\btheta}}

\def\RE{\operatorname{RE}}
\def\bSigma{\bolds{\Sigma}}
\def\Sigmabar{{\bar{\Sigma}}}
\def\bSigmabar{{\bar{\bSigma}}}
\def\hbG{{\hat\bG}}
\def\bmu{{\bolds\mu}}
\def\bDelta{{\bolds\Delta}}

\makeatother

\begin{document}
\begin{frontmatter}

\title{Oracle inequalities for the lasso in the Cox model}
\runtitle{Lasso in the Cox model}

\begin{aug}
\author[A]{\fnms{Jian} \snm{Huang}\thanksref{t1}\ead[label=e1]{jian-huang@uiowa.edu}},
\author[B]{\fnms{Tingni} \snm{Sun}\ead[label=e2]{tingni@wharton.upenn.edu}},
\author[C]{\fnms{Zhiliang} \snm{Ying}\thanksref{t3}\ead[label=e3]{zying@stat.columbia.edu}},
\author[D]{\fnms{Yi} \snm{Yu}\ead[label=e4]{yuyi@fudan.edu.cn}}\\
\and
\author[E]{\fnms{Cun-Hui} \snm{Zhang}\corref{}\thanksref{t5}\ead[label=e5]{czhang@stat.rutgers.edu}\ead[label=u1,url]{http://www.foo.com}}
\runauthor{J. Huang et al.}
\affiliation{University of Iowa, University of Pennsylvania,
Columbia University, Fudan University and Rutgers University}
\address[A]{J. Huang\\
Department of Statistics\\
\quad and Actuarial Science\\
University of Iowa \\
241 SH \\
Iowa City, Iowa 52242\\
USA\\
\printead{e1}}
\address[B]{T. Sun\\
Statistics Department\\
The Wharton School\\
University of Pennsylvania\\
400 Jon M. Huntsman Hall\\
3730 Walnut Street \\
Philadelphia, Pennsylvania 19104-6340\\
USA\\
\printead{e2}}
\address[C]{Z. Ying\\
Department of Statistics\\
Columbia University\\
1255 Amsterdam Avenue\\
New York, New York 10027\\
USA\\
\printead{e3}}
\address[D]{Y. Yu\\
School of Mathematical Sciences\hspace*{19.65pt}\\
Fudan University\\
Shanghai\\
China\\
\printead{e4}}
\address[E]{C.-H. Zhang \\
Department of Statistics\\
\quad and Biostatistics \\
Hill Center, Busch Campus\\
Rutgers University\\
Piscataway, New Jersey 08854\\
USA\\
\printead{e5}} 
\end{aug}

\thankstext{t1}{Supported in part by NIH Grants R01CA120988,
R01CA142774 and NSF Grants DMS-08-05670 and DMS-12-08225.}
\thankstext{t3}{Supported in part by NIH Grant R35GM047845 and NSF
Grant SES1123698.}
\thankstext{t5}{Supported in part by NSF Grants DMS-09-06420, DMS-11-06753
and DMS-12-09014 and NSA Grant H98230-11-1-0205.}

\received{\smonth{11} \syear{2012}}
\revised{\smonth{2} \syear{2013}}

%
\begin{abstract}
We study the absolute penalized maximum partial likelihood estimator in
sparse, high-dimensional Cox proportional hazards regression models
where the number of time-dependent covariates can be larger than the
sample size. We establish oracle inequalities based on natural
extensions of the compatibility and cone invertibility factors of the
Hessian matrix at the true regression coefficients. Similar results
based on an extension of the restricted eigenvalue can be also proved
by our method. However, the presented oracle inequalities are sharper
since the compatibility and cone invertibility factors are always
greater than the corresponding restricted eigenvalue. In the Cox
regression model, the Hessian matrix is based on time-dependent
covariates in censored risk sets, so that the compatibility and cone
invertibility factors, and the restricted eigenvalue as well, are
random variables even when they are evaluated for the Hessian at the
true regression coefficients. Under mild conditions, we prove that
these quantities are bounded from below by positive constants for
time-dependent covariates, including cases where the number of
covariates is of greater order than the sample size. Consequently, the
compatibility and cone invertibility factors can be treated as positive
constants in our oracle inequalities.
\end{abstract}

%

\end{frontmatter}

\section{Introduction}\label{sec1}
The \citet{Cox72} proportional hazards model is widely used in the
regression analysis of censored
survival data, notably in identifying risk factors in epidemiological
studies and treatment effects in
clinical trials when the outcome variable is time to event. In a
traditional biomedical study,
the number of covariates $p$ is usually relatively small as compared
with the sample size $n$.
Theoretical properties of the maximum partial likelihood estimator in
the fixed $p$ and large $n$
setting have been well established.
For example, \citet{Tsi81} proved the asymptotic normality of the
maximum partial likelihood estimator.
\citet{AndGil82} formulated the Cox model in the context of the
more general counting process
framework and studied the asymptotic properties of the maximum partial
likelihood estimator
using martingale techniques. These results provide a solid foundation
for the applications of the Cox model in a diverse range of problems
where time to event is the outcome of interest.

In recent years, technological advancement has resulted in the
proliferation of
massive high-throughput and high-dimensional
genomic data in studies that attempt to find genetic risk factors for
disease and clinical outcomes,
such as the age of disease onset or time to death.
Finding genetic risk factors for survival is fundamental to modern
biostatistics,
since survival is an important clinical endpoint.
However, in such problems, the standard approach to the Cox model is
not applicable,
since the number of potential genetic risk factors is typically
much larger than the sample size.
In addition, traditional variable selection methods such as subset
selection are not computationally
feasible when $p \gg n$.

The $\ell_1$-penalized least squares estimator, or the lasso, was
introduced by \citet{Tib96}.
In the wavelet setting, the $\ell_1$-penalized method was introduced
by \citet{CheDonSau98}
as basis pursuit.
Since then, the lasso has emerged as a widely used approach to variable
selection and estimation in sparse,
high-dimensional statistical problems.
It has also been extended to the Cox model [\citet{Tib97}].
\citet{GuiLi05} implemented the LARS algorithm [\citet{Efretal04}]
to approximate the lasso
in the Cox regression model and applied their method to survival data
with microarray gene expression covariates.
Their work demonstrated the effectiveness of the lasso for variable
selection in the Cox model in a $p \gg n$
setting.

There exists a substantial literature on the lasso and other penalized
methods for survival models with a fixed number of covariates $p$.
\citet{ZhaLu07} considered an adaptive lasso for the Cox model and
showed that, under certain
regularity conditions and with a suitable choice of the penalty
parameter, their method possesses
the asymptotic oracle property when the maximum partial likelihood
estimator is used as the initial estimator.
\citet{FanLi02} proposed the use of the smoothly clipped absolute
deviation (SCAD) penalty
[\citet{Fan97}, \citet{FanLi01}] for variable selection and estimation
in the Cox model which may include a frailty term.
With a suitable choice of the penalty parameter, they showed that a
local maximizer of the penalized log-partial
likelihood has an asymptotic oracle property
under certain regularity conditions on the Hessian of the log-partial
likelihood and the censoring mechanism.



Extensive efforts have been focused upon the analysis of regularization
methods for
variable selection in the $p \gg n$ setting.
In particular, considerable progress has been made in theoretical
understanding of the lasso.
However, most results are developed in the linear regression model.
\citet{GreRit04} studied the prediction performance of the
lasso in high-dimensional
linear regression.
\citet{MeiBuh06} showed that, for neighborhood
selection in the Gaussian graphical model,
under a neighborhood stability condition and certain additional
regularity conditions,
the lasso is consistent even when $p/n \rightarrow\infty$.
\citet{ZhaYu06} formalized the neighborhood stability condition in
the context of linear regression
as a strong irrepresentable condition on the design matrix.
Oracle inequalities for the prediction and estimation error of the
lasso was developed in
\citet{BunTsyWeg07}, \citet{ZhaHua08}, \citet{MeiYu09},
\citet{BicRitTsy09}, \citet{Zha09} and \citet{YeZha10},
among many others.

A number of papers analyzed penalized methods beyond linear regression.
\citet{FanPen04} established oracle properties for a local
solution of
concave penalized estimator in a general setting with $n\gg p\to\infty$.
\citet{van08} studied the lasso in high-dimensional generalized
linear models (GLM)
and obtained prediction and $\ell_1$ estimation error bounds.
\citet{Negetal10} studied penalized
M-estimators with a general class of
regularizers, including an $\ell_2$ error bound for the lasso in GLM
under a restricted convexity and other
regularity conditions.
\citet{BraFanJia11} made significant progress by extending the results
of \citet{FanLi01}
to a more general class of penalties in the Cox regression model with
large $p$
under different sets of regularity conditions.
\citet{HuaZha12} studied weighted absolute penalty
and its adaptive, multistage application in GLM.


In view of the central role of the Cox model in survival analysis, its
widespread applications and
the proliferation of $p\gg n$ data, it is of great interest to
understand the properties of the related lasso approach.
The main goal of the present paper is to establish theoretical
properties for the lasso
in the Cox model when \mbox{$p \gg n$}.
Specifically, we extend certain basic inequalities from linear
regression to the case of the Cox regression.
We generalize the compatibility and cone invertibility factors from the
linear regression
model and establish oracle inequalities for the lasso in the Cox
regression model
in terms of these factors at the true parameter value.
Moreover, we prove that the compatibility and cone invertibility
factors can be treated
as constants under mild regularity conditions.

A main feature of our results is that they are derived under the more
general counting process formulation
of the Cox model with potentially a larger number of time-dependent
covariates than the sample size.
This formulation is useful because it ``permits a regression analysis
of the intensity of a recurrent event allowing
for complicated censoring patterns and time-dependent covariates''
[\citet{AndGil82}].

A second main feature of our results is that the regularity conditions
on the counting
processes and time-dependent covariates are directly imposed on the
compatibility and cone invertibility
factors of the Hessian of the negative log-partial likelihood evaluated
at the true regression coefficients.
Under such regularity conditions, the lasso estimator is proven to live
in a neighborhood where
the ratio between the estimated and true hazards is uniformly bounded
away from zero
and infinity.
This allows unbounded and near zero ratios between the true and
baseline hazards.
Our analysis can be also used to prove oracle inequalities based on the
restricted eigenvalue.
However, since the compatibility and cone invertibility factors are
greater than the corresponding
restricted eigenvalue [\citet{vanBuh09}, \citet{YeZha10}],
the presented results
are sharper.

A third main feature of our results is that the compatibility and cone
invertibility factors used, and
the smaller corresponding restricted eigenvalue, are proven to be
greater than a fixed positive
constant under mild conditions on the
counting processes and time-dependent covariates,
including cases where $p\gg n$.
In the Cox regression model, the Hessian matrix is based on weighted
averages of the cross-products
of time-dependent covariates in censored risk sets,
so that the compatibility and cone invertibility factors and the
restricted eigenvalue are random
variables even when they are evaluated for the Hessian at the true
regression coefficients.
Under mild conditions, we prove that these quantities are bounded from
below by
positive constants as certain truncated population versions of them.
Thus the compatibility and cone invertibility factors can be treated as
constants in our oracle inequalities.

The main results of this paper and the analytical methods used
for deriving them are identical to those in its
predecessor submitted in November 2011, with Section~\ref{sec4} as an exception.
The difference in Section~\ref{sec4} is that the lower bound for the
compatibility and cone invertibility factors and the restricted
eigenvalue is improved
to allow time-dependent covariates.

During the revision process of our paper,
we became aware of a number of papers on hazard regression
with censored data.
\citet{KonNan} took an approach of \citet{van08} to derive
prediction and
$\ell_1$ error bounds for the lasso in the Cox proportional hazards regression
under a quite different set of conditions from us. For example, they
required an $\ell_1$ bound
on the regression coefficients to guarantee a uniformly bounded ratio
between hazard functions
under consideration.
\citet{Lem} considered the joint estimation of the baseline hazard
function and regression
coefficients in the Cox model.
As a result, Lemler's (\citeyear{Lem}) error bounds for regression coefficients
are of greater order than ours
when the intrinsic dimension of the unknown baseline hazard function
is of greater order than the number of nonzero regression coefficients.
\citet{GafGui12} considered a quadratic loss function in
place of a negative log-likelihood function
in an additive hazards model. A nice feature of the additive hazards
model is that the quadratic loss
actually produces unbiased linear estimation equations so that the
analysis of the lasso is
similar to that of linear regression.
The oracle inequalities in these three papers and ours can be all
viewed as nonasymptotic.
Unlike our paper, none of these three papers consider time-dependent
covariates or
constant lower bounds of the restricted eigenvalue or related key factors
for the analysis of the lasso.

The rest of this paper is organized as follows.
In Section~\ref{sec2} we provide basic notation and model specifications.
In Section~\ref{sec3} we develop oracle inequalities for the lasso in
the Cox model.
In Section~\ref{sec4} we study the compatibility and cone invertibility factors
and the corresponding restricted eigenvalue
of the Hessian of the log-partial likelihood in the Cox model.
In Section~\ref{sec5} we make some additional remarks.
All proofs are provided either right after the statement of the result
or deferred to the \hyperref[app]{Appendix}.

\section{\texorpdfstring{Cox model with the $\ell_1$ penalty}
{Cox model with the l1 penalty}}\label{sec2}

Following \citet{AndGil82}, consider an $n$-dimensional
counting process
$\mathbf{N}^{(n)} (t)= (N_1(t),\ldots, N_n(t))$, $t\ge0$,
where $N_i(t)$ counts the number of observed events for the $i$th
individual in the time interval $[0, t]$.
The sample paths of $N_1,\ldots, N_n$ are step functions, zero at
$t=0$, with jumps of size $+1$ only.
Furthermore, no two components jump at the same time.
For $ t\ge0$, let ${\cal F}_t$ be the $\sigma$-filtration
representing all the information available up to time $t$.
Assume that for $\{{\cal F}_t, t\ge0\}$, $\mathbf{N}^{(n)}$ has
predictable compensator $\bolds{\Lambda}^{(n)} = (\Lambda_1,\ldots,
\Lambda_n)$ with
%
%
\begin{equation}
\label{ph-model} d\Lambda_i(t) = Y_i(t)\exp\bigl\{
\bZ'_i(t){\bbeta}^o\bigr\}\,d
\Lambda_0(t),
\end{equation}
where $\bbeta^o$ is a $p$-vector of true regression coefficients,
$\Lambda_0$ is an unknown baseline cumulative hazard function and, for
each $i$,
$Y_i(t)\in\{0, 1\}$ is a predictable at risk indicator process that
can be constructed from data,
and $\mathbf{Z}_i(t)=(Z_{i,1}(t),\ldots, Z_{i,p}(t))'$
is a $p$-dimensional vector-valued predictable covariate process.
In this setting the $\sigma$-filtration can be the natural
${\cal F}_t=\sigma\{N_i(s), Y_i(s),\bZ_i(s);\allowbreak s \le t, i=1,\ldots,
n\}$ or a richer one.
We are interested in the problem of variable selection in sparse,
high-dimensional settings where~$p$,
the number of possible covariates, is large, but the number of
important covariates is relatively small.

\subsection{\texorpdfstring{Maximum partial likelihood estimator with $\ell_1$ penalty}
{Maximum partial likelihood estimator with l1 penalty}}\label{sec2.1}
Define logarithm of the Cox partial likelihood for survival experience
at time $t$,
\[
C(\bbeta;t) = \sum_{i=1}^n \int
_0^t \bZ'_i(s){
\bbeta} \,dN_i(s) - \int_0^t \log
\Biggl[\sum_{i=1}^n Y_i(s)
e^{\bZ'_i(s){\bbeta
}} \Biggr] \,d\Nbar(s),
\]
where $\Nbar=\sum_{i=1}^n N_i$. The log-partial likelihood function
is
\[
C(\bbeta,\infty)=\lim_{t\to\infty}C(\bbeta,t).
\]
Let $\ell(\bbeta) = -C(\bbeta;\infty)/n$.
The maximum partial likelihood estimator is the value that minimizes
$\ell(\bbeta)$.

An approach to variable selection in sparse, high-dimensional settings
for the Cox model is to minimize an $\ell_1$-penalized negative
log-partial likelihood criterion,
%
%
\begin{equation}
\label{pen-lik} L(\bbeta; \lambda) = \ell(\bbeta) + \lam|\bbeta|_1
\end{equation}
[\citet{Tib97}], where $\lam\ge0$ is a penalty parameter.
Henceforth, we use notation $|\bbeta|_q= \{\sum_{i=1}^p|\beta_i|^q \}
^{1/q}$ for $1\le q< \infty$,
$|\bbeta|_{\infty}=\max_{1\le i\le p}|\beta_i|$ and $|\bbeta|_0=\#
\{j\dvtx\beta_j\neq0\}$. For a given $\lam$, the $\ell_1$-penalized
maximum partial likelihood estimator, or the lasso estimator hereafter,
is defined as
%
%
\begin{equation}
\label{ph-Lasso} \bhbeta(\lam) = \argmin_{\bbeta} L(\bbeta;\lam).
\end{equation}

\subsection{The Karush--Kuhn--Tucker conditions}\label{sec2.2}
The lasso estimator can be characterized by the Karush--Kuhn--Tucker
(KKT) conditions. Since the log-partial likelihood belongs to an
exponential family, $\ell(\bbeta)$ must be convex in $\bbeta$ and so is
$L(\bbeta;\lam)$. It follows that a vector
$\bhbeta=(\hbeta_1,\ldots,\hbeta_p)'$ is a solution to (\ref{ph-Lasso})
if and only if the following KKT conditions hold:
%
%
\begin{equation}
\label{KKT} \cases{\dell_j(\bhbeta) =- \lam\sgn(
\hbeta_j), &\quad if $\hbeta_j \neq0$,
\vspace*{2pt}\cr
\bigl|
\dell_j(\bhbeta)\bigr| \le\lam, &\quad if $\hbeta_j=0$,}
\end{equation}
where $\dell(\bbeta)=(\dell_1(\bbeta),\ldots,\dell_p(\bbeta))'
= \partial\ell(\bbeta)/\partial\bbeta$
is the gradient of $\ell(\bbeta)$.
The necessity and sufficiency of (\ref{KKT}) can be proved by
subdifferentiation of the convex penalized loss (\ref{pen-lik}).
This does not require strict convexity.

The KKT conditions indicate that the lasso in the Cox regression model
may be analyzed in a similar way to the
lasso in linear regression. As can be seen in the subsequent
developments, such analysis can be carried out by proving that
$|\dell(\bbeta^o)|_\infty$ is sufficiently small
and the Hessian of $\ell(\bbeta)$ does not vanish for
a sparse
$\bbeta$ at the true $\bbeta= {\bbeta}^o$.
The (local) martingales for the counting process will play a crucial
role to ensure that
these requirements are satisfied.

\subsection{Additional notation}\label{sec2.3}
Since the $\Lambda_i$ are compensators,
\[
M_i(t)=N_i(t)-\int_0^tY_i(s)
\exp\bigl(\bZ'_i(s){\bbeta}^o\bigr) \,d
\Lambda_0(s),\qquad 1\le i \le n, t\ge0,
\]
are (local) martingales with predictable variation/covariation processes
\[
\langle M_i, M_i\rangle(t) =\int_0^t
Y_i(s)\exp\bigl(\bZ'_i(s){\bbeta
}^o\bigr) \,d\Lambda_0(s) \quad\mbox{and}\quad\langle
M_i, M_j\rangle= 0,\qquad i \neq j.\vadjust{\goodbreak}
\]
For a vector $v$, let $v^{\otimes0}=1\in\real$, $v^{\otimes1}=v$
and $v^{\otimes2}=vv'$. Define
\begin{eqnarray*}
S^{(k)}(t,\bbeta) & = & \frac{1}{n}\sum
_{i=1}^n \mathbf{Z}_i^{\otimes k}(t)
Y_i(t) e^{\bZ'_i(t){\bbeta}},\qquad k=0,1,2,
\\
R_n(t,\bbeta) &=& \frac{1}{n}\sum
_{i=1}^n Y_i(t) e^{\bZ'_i(t){\bbeta}},\qquad
\bar{
\mathbf{Z}}_n(t,\bbeta)=\frac{S^{(1)}(t,
\bbeta)}{S^{(0)}(t,\bbeta)},
\\
V_n(t,{\bbeta}) &=& \sum
_{i=1}^n w_{ni}(t,{\bbeta}) \bigl(\mathbf
{Z}_i(t)-\bar{\mathbf{Z}}_n(t,\bbeta)
\bigr)^{\otimes2} =\frac{S^{(2)}(t,\bbeta)}{S^{(0)}(t,\bbeta)} - \bar
{\mathbf{Z}}_n(t,{
\bbeta})^{\otimes2},
\end{eqnarray*}
where $w_{ni}(t,{\bbeta})=Y_i(t) \exp(\bZ'_i(t){\bbeta
})/[nS^{(0)}(t,{\bbeta})]$.
By differentiation and rearrangement of terms, it can be shown as in
\citet{AndGil82} that the gradient of
$\ell({\bbeta})$ is
%
%
\begin{equation}
\label{gradient1} \dell({\bbeta}) \equiv
\frac{\partial\ell({\bbeta})}{\partial{\bbeta}} = -
\frac{1}{n}\sum_{i=1}^n \int
_0^\infty\bigl[\mathbf{Z}_i(s)- \bar{
\mathbf{Z}}_n(s,{\bbeta})\bigr]\,dN_i(s), 
\end{equation}
and the Hessian matrix of $\ell({\bbeta})$ is
%
%
\begin{equation}
\label{hessian1} 
\ddell({\bbeta}) \equiv\frac{\partial^2\ell({\bbeta})}{\partial
{\bbeta}\,\partial{\bbeta}'} =
\frac{1}{n}\int_0^\infty V_n(s,{
\bbeta}) \,d\Nbar(s).
\end{equation}

\section{Oracle inequalities}\label{sec3}
In this section, we derive oracle inequalities for the estimation error
of lasso in the Cox regression model.
Let ${\bbeta}^o$ be the vector of true regression coefficients, and define
$\cO=\{j\dvtx \beta_j^o\neq0\}$, $\cO^c=\{j\dvtx\break \beta_j^o= 0\}$ and
$d_o=|\cO|$,
where $|\mathcal{U}|$ for a set $\mathcal{U}$ denotes its cardinality.

Making use of the KKT conditions (\ref{KKT}), we first develop a basic
inequality involving
the symmetric Bregman divergence and $\ell_1$ estimation error in the
support $\cO$ of $\bbeta^o$ and its
complement. The symmetric Bregman divergence, defined as
\[
D^s(\bhbeta,\bbeta) = (\bhbeta- \bbeta)' \bigl(
\dell(\bhbeta)-\dell(\bbeta) \bigr)
\]
can be viewed as symmetric, partial Kullback--Leibler distance between
the partial likelihood at
$\bhbeta$ and $\bbeta$.
Thus, $D^s(\bhbeta,\bbeta)$ can be viewed as a measure
of prediction performance.
The basic inequality, given in Lemma~\ref{LemH1} below, serves as a vehicle
for establishing the desired oracle inequalities.


%
\begin{lemma}
\label{LemH1}
Let $\bhbeta$ be defined as in (\ref{ph-Lasso}),
$\tbtheta=\bhbeta-{\bbeta}^o$ and $z^* =
|\dell(\bbeta^o)|_{\infty}$. Then the following
inequalities hold:
%
%
\begin{equation}
\label{LemH1E1} \bigl(\lam- z^*\bigr) |\tbtheta_{\cO^c}|_1
\le D^s(\bhbeta,\bbeta) 
+ \bigl(\lam- z^*
\bigr) |\tbtheta_{\cO^c}|_1 \le\bigl(\lam+ z^*\bigr) |
\tbtheta_{\cO}|_1,
\end{equation}
where $\tbtheta_{\cO}$ and $\tbtheta_{\cO
^c}$ denote the subvectors of $\tbtheta$ of components in
$\cO$ and $\cO^c$, respectively.
In particular, for any $\xi>1$,
$|\tbtheta_{\cO^c}|_1\le\xi|\tbtheta_{\cO}|_1$
in the event $z^*\le(\xi-1)/(\xi+1)\lam$.
\end{lemma}

It follows from Lemma~\ref{LemH1} that in the event $z^*\le(\xi
-1)/(\xi+1)\lam$,
the estimation error $\tbtheta=\bhbeta
-\bbeta^o$ belongs to the cone
%
%
\begin{equation}
\label{cone} \mathscr C(\xi,\cO) = \bigl\{\mathbf{b} \in\real^p\dvtx |
\mathbf{b}_{\cO^c}|_1 \le\xi|\mathbf{b}_{\cO}|_1
\bigr\}.
\end{equation}
In linear regression, the invertibility of the Gram matrix in the same
cone, expressed in terms of
restricted eigenvalues and related quantities, has been used to control
the estimation error of the lasso.
In what follows, we prove that a direct extension of the compatibility
and cone invertibility factors
can be used to control the estimation error of the lasso in the Cox regression.

For the cone in (\ref{cone}) and a given $p\times p$
nonnegative-definite matrix $\bSigmabar$, define
%
%
\begin{equation}
\label{Re1} \kappa(\xi,\cO;\bSigmabar) = \inf_{0\neq\mathbf{b}\in
\mathscr C(\xi, \cO)}
\frac
{d_o^{1/2}(\mathbf{b}'\bSigmabar\mathbf{b})^{1/2}} {
|\mathbf{b}_{\cO}|_1}
\end{equation}
as the compatibility factor [\citet{van07}, \citet
{vanBuh09}], and
%
%
\begin{equation}
\label{Fq} F_q(\xi,\cO;\bSigmabar) = \inf_{0\neq\mathbf{b}\in\mathscr
C(\xi, \cO)}
\frac
{d_o^{1/q}\mathbf{b}'\bSigmabar\mathbf{b}} {
|\mathbf{b}_{\cO}|_1|\mathbf{b}|_q}
\end{equation}
as the weak cone invertibility factor [\citet{YeZha10}].
These quantities are closely related to the restricted eigenvalue
[\citet{BicRitTsy09}, \citet{Kol09}],
%
%
\begin{equation}
\label{Re2} \RE(\xi,\cO;\bSigmabar) =\inf_{0\neq\mathbf{b}\in\mathscr
C(\xi, \cO)}
\frac
{(\mathbf{b}'\bSigmabar\mathbf{b})^{1/2}} {
|\mathbf{b}|_2}.
\end{equation}

In linear regression, the Hessian of the squared loss $|\by-\bX
\bbeta|_2^2/(2n)$ is
taken as~$\bSigmabar$, and the oracle inequalities established in the
papers cited in the above
paragraph can be summarized as follows:
in the event $z^* = |\bX'(\by-\bX\bbeta^o)/n|_\infty\le\lam(\xi
-1)/(\xi+1)$,
%
%
\begin{equation}
\label{LM1} \bigl|\bX\bigl(\hbbeta-{\bbeta}^o\bigr)\bigr|_2^2/n
\le\frac{4(1+1/\xi)^{-2}\lam^2 d_o}{\kappa^2(\xi,\cO;\bX'\bX/n)},\qquad
\bigl|\hbbeta-{\bbeta}^o\bigr|_1\le\frac{2\xi
d_o\lam}{\kappa^2(\xi,\cO;\bX'\bX/n)}\hspace*{-32pt}
\end{equation}
and
%
%
\begin{eqnarray}
\label{LM2}
\bigl|\hbbeta-{\bbeta}^o\bigr|_2&\le&\frac{2(1+1/\xi)^{-1}d_o^{1/2}\lam
}{\RE
^2(\xi,\cO;\bX'\bX/n)},\nonumber\\[-8pt]\\[-8pt]
\bigl|\hbbeta-{\bbeta}^o\bigr|_q&\le&\frac{2(1+1/\xi)^{-1}d_o^{1/q}\lam
}{F_q(\xi,\cO;\bX'\bX/n)},\qquad q\ge1.\nonumber
\end{eqnarray}

In the Cox regression model, we still take the Hessian of the
log-partial likelihood as $\bSigmabar$,
in fact the Hessian at the true $\bbeta^o$,
so that (\ref{Re1}) and (\ref{Fq}) become
%
%
\begin{equation}
\label{Fq-Cox} \kappa(\xi,\cO) = \kappa\bigl(\xi,\cO;\ddell\bigl(
\bbeta^o\bigr)\bigr),\qquad F_q(\xi,\cO) = F_q
\bigl(\xi,\cO;\ddell\bigl(\bbeta^o\bigr)\bigr).
\end{equation}
The reason for using these factors is that they yield somewhat sharper
oracle inequalities than
the restricted eigenvalue.
It follows from $|\bb_{\cO}|_1\le d_o^{1/2}|\bb|_2$ that
$F_2(\xi,\cO)\ge\kappa(\xi,\cO)\RE(\xi,\cO)$ and $\kappa(\xi,\cO)\ge\RE
(\xi,\cO)$.
Therefore, the first inequality of (\ref{LM2}) is subsumed by the
second with $q=2$.
Moreover, it is possible to have $\kappa(\xi,\cO)\gg\RE(\xi,\cO)$
[\citet{vanBuh09}], and consequently, the $\ell_2$
error bound
based on the cone invertibility factor may be of sharper order that the
one based on the restricted eigenvalue.

The following theorem extends (\ref{LM1}) and (\ref{LM2})
from the linear regression model to the proportional hazards regression
model with
%
%
\begin{equation}
\label{cond-K} \max_{i < i' \le n}\sup_{0\le t<\infty}\max
_{j\le p}\bigl|Z_{i,j}(t) - Z_{i',j}(t)\bigr|\le K.
\end{equation}
Let $\xi>1$, $\cO=\{j\dvtx \beta^o_j\neq0\}$,
$\kappa(\xi,\cO)$ and $F_q(\xi,\cO)$ be as in (\ref{Fq-Cox}).

%
\begin{theorem}
\label{th-1} Let $\tau=K(\xi+1) d_o\lambda/\{2\kappa^2(\xi,\cO)\}
$ with a certain $K>0$.
Suppose condition (\ref{cond-K}) holds and $\tau\le1/e$.
Then, in the event $|\dell({\bbeta}^o)|_{\infty}\le(\xi-1)/\allowbreak(\xi
+1)\lam$,
%
%
\begin{equation}
\label{th-1-1}\quad D^s(\bhbeta,\bbeta) \le\frac{4e^\eta(1+1/\xi
)^{-2}\lam^2 d_o}{\kappa^2(\xi,\cO)},\qquad 
\bigl|\bhbeta-{\bbeta}^o\bigr|_1
\le\frac{e^{\eta}(\xi+1)d_o\lam}{2\kappa^2(\xi,\cO)}
\end{equation}
and
%
%
\begin{equation}
\label{th-1-2} \bigl|\bhbeta-{\bbeta}^o\bigr|_q\le
\frac{e^\eta2(1+1/\xi)^{-1} d_o^{1/q}\lam}{F_q(\xi,\cO)},\qquad q\ge1,
\end{equation}
where $\eta\le1$ is the smaller solution of $\eta e^{-\eta}=\tau$.
\end{theorem}

Compared with (\ref{LM1}) and (\ref{LM2}), the new inequalities (\ref
{th-1-1}) and (\ref{th-1-2})
contain an extra factor $e^\eta\le e$. This is due to the nonlinearity
in the Cox regression score equation.
Aside from this factor, the error bounds for the Cox regression have
the same form as those for linear
regression, except for an improvement of a factor of $4\xi/(1+\xi)\ge
2$ for the $\ell_1$ oracle inequality.

The theorem assumes condition (\ref{cond-K}), which asserts $|\bZ
_i(t)-\bZ_{i'}(t)|_\infty\le K$ uniformly in $\{t,i,i'\}$.
This condition is a consequence of the uniform boundedness of the
individual covariates,
and is reasonable in most practical situations (e.g., single-nucleotide
polymorphism data).
In the case where the covariates are normal variables with uniformly
bounded variance,
the condition holds with $K=K_{n,p}$ of $\sqrt{\log(np)}$ order.

From an analytical perspective, an important feature of (\ref{th-1-1})
and (\ref{th-1-2}) is that the constant factors
(\ref{Re1}) and (\ref{Fq}) are both defined with the true ${\bbeta
}^o$ in (\ref{Fq-Cox}).
No condition is imposed on the gradient and Hessian of the log-partial
likelihood for ${\bbeta}\neq{\bbeta}^o$.
In other words, the key condition $\tau<1/e$, expressed in terms of $\{
K,d^o,\lam\}$ and the
compatibility factor $\kappa^2(\xi,\cO)$ at the true $\bbeta^o$, is
sufficient to guarantee\vadjust{\goodbreak}
the error bounds in Theorem~\ref{th-1}.
Thus, our results are much simpler to state and conditions easier to
verify than existing ones requiring
regularity conditions in a neighborhood of $\bbeta^o$ in the Cox
regression model.
This feature of Theorem~\ref{th-1} plays a crucial role in our
derivation of lower bounds for
$\kappa^2(\xi,\cO)$ and $F_q(\xi,\cO)$ for time-dependent
covariates in Section~\ref{sec4}.
We note that the local martingale structure is valid only at the true
${\bbeta}^o$.

To prove Theorem~\ref{th-1}, we develop a sharpened version of an
inequality of \citet{HjoPol93}.
This inequality, given in Lemma~\ref{LemB} below, explicitly controls
the symmetric Bregman-divergence and
Hessian of the log-partial likelihood in a neighborhood of $\bbeta$.
Based on this relationship, Theorem~\ref{th-1} is proved using the
definition of the
quantities in (\ref{Fq-Cox}) and the membership of the error
$\bhbeta-{\bbeta}^o$ in the cone $\mathscr C(\xi,\cO)$
(\ref{cone}).
For two symmetric matrices $A$ and $B$, $A \le B$ means $B-A$ is
nonnegative-definite.

%
\begin{lemma}
\label{LemB}
Let $\ell({\bbeta})$ and its Hessian $\ddell({\bbeta})$ be as in
(\ref{pen-lik})
and (\ref{hessian1}). Then
%
%
\begin{equation}
\label{LemB1}\qquad e^{-\eta_{\mathbf{b}}}\mathbf{b}'\ddell({\bbeta})
\mathbf{b} \le D^s(\bbeta+\bb,\bbeta) =\mathbf{b}'\bigl[
\dell({\bbeta}+\mathbf{b}) - \dell({\bbeta})\bigr] \le e^{\eta_{\mathbf{b}}}
\mathbf{b}'\ddell({\bbeta})\mathbf{b},
\end{equation}
where $\eta_{\mathbf{b}}=\max_{s\ge0}\max_{i,j}|\mathbf{b}'
\mathbf{Z}_i(s)-\mathbf{b}'\mathbf{Z}_j(s)|$. Moreover,
%
%
\begin{equation}
\label{LemB2} e^{-2\eta_{\mathbf{b}}} \ddell({\bbeta}) \le\ddell({\bbeta
}+\mathbf{b})
\le e^{2\eta_{\mathbf
{b}}}\ddell({\bbeta}).
\end{equation}
\end{lemma}

Under the conditions of Theorem~\ref{th-1}, the factors $e^{\pm\eta
_{\mathbf{b}}}$ and
$e^{\pm2 \eta_{\mathbf{b}}}$ in the inequalities in Lemma~\ref{LemB}
are bounded by positive constants.
These factors lead to the factor $e^{\eta}$ for $\eta\le1$ (and thus
$e^{\eta} \le e$) in the upper bounds in (\ref{th-1-1}) and (\ref{th-1-2}).

Since the oracle inequalities in Theorem~\ref{th-1} are guaranteed to
hold only
within the event $|\dell({\bbeta}^o)|_{\infty}\le(\xi-1)/(\xi
+1)\lam$,
a probabilistic upper bound is needed for $|\dell({\bbeta
}^o)|_{\infty}$.
Lemma~\ref{LemC} below provides such a probability bound.
Similar inequalities can be found in \citet{del99}.

%
\begin{lemma}
\label{LemC}
\textup{(i)} Let $f_n(t) = n^{-1}\sum_{i=1}^n \int_0^t a_i(s)\{dN_i(t)-
Y_i(s)\exp(\bZ_i'(s)\*\bbeta^o)\,d\Lambda_0(s)\}$
with $[-1,1]$-valued predictable processes $a_i(s)$. Then,
for all\break \mbox{$C_0>0$},
%
%
\begin{eqnarray}\quad
\label{lm-C-0} \rP\Biggl\{\max_{t>0}\bigl|f_n(t)\bigr| &>&
C_0 x, \sum_{i=1}^n \int
_0^\infty Y_i(t) \,dN_i(t)
\le C_0^2n \Biggr\} \le2 e^{- n x^2/2}.
\end{eqnarray}

\textup{(ii)} Suppose that $\max_{i\le n}\sup_{t\ge0}\max_{j\le
p}|Z_{i,j}(t)-\bar{Z}_{n,j}(t,{\bbeta}^o)|_\infty\le K$,
where $\bar{Z}_{n,j}(t,{\bbeta}^o)$ are the components of ${\bar\bZ
}_{n}(t,{\bbeta}^o)$.
Let $\dell({\bbeta})$ be the gradient in (\ref{gradient1}). Then,
for all $C_0>0$,
%
%
\begin{eqnarray}
\label{lm-C-1}\quad \rP\Biggl\{\bigl|\dell\bigl({\bbeta}^o
\bigr)\bigr|_\infty> C_0Kx, \sum_{i=1}^n
\int_0^\infty Y_i(t)
\,dN_i(t) \le C_0^2n \Biggr\} &\le&2 p
e^{- n x^2/2}.
\end{eqnarray}
In particular, if $\max_{i\le n}N_i(1)\le1$, then
$\rP\{|\dell({\bbeta}^o)|_\infty> Kx\}\le2 p e^{- n x^2/2}$.\vadjust{\goodbreak}
\end{lemma}

The following theorem states an upper bound of the estimation error,
which follows directly from Theorem~\ref{th-1} and Lemma~\ref{LemC}.

%
\begin{theorem}
\label{ThmA}
Suppose (\ref{cond-K}) holds and $N_i(\infty)\le1$ for all
$i\le n$ and $t\ge0$. Let $\xi>1$ and $\lambda= \{(\xi+1)/(\xi-1)\}
K\sqrt{(2/n)\log(2p/\eps)}$
with a small $\varepsilon>0$ (e.g., $\eps=1\%$). Let $C_{\kappa}>0$ satisfying
$\tau= K(\xi+1)\,d_o\lambda/(2C_{\kappa}^2) \le1/e$. Let
$\eta\le1$ be the smaller solution of $\eta e^{-\eta}=\tau$.
Then, for any $C_{F,q}>0$,
\begin{eqnarray*}
D^s(\bhbeta,\bbeta) &\le&\frac{4e^\eta\xi^2\lam^2
d_o}{(1+\xi)^2C_{\kappa}^2},\qquad \bigl|\bhbeta-{
\bbeta}^o\bigr|_1\le\frac{e^{\eta}(\xi
+1)d_o\lambda}{2C_{\kappa}^2},\\
\bigl|\bhbeta-
\bbeta^o\bigr|_q&\le&\frac{2e^\eta\xi
d_o^{1/q}\lambda}{(\xi+1)C_{F,q}}
\end{eqnarray*}
all hold with at least probability $\rP\{\kappa(\xi,\cO) \ge
C_{\kappa}, F_q(\xi,\cO) \ge C_{F,q}\} - \eps$.
\end{theorem}

It is noteworthy that this theorem gives an upper bound of the
estimation error for all the $\ell_q$ norms with $q \ge1$.
From this theorem, for the $\ell_q$ error $|\bhbeta
-{\bbeta}^o|_q$
with $q \ge1$ to be small with high probability, we need to ensure that
$d_o \lambda\to0$ as $n\to\infty$.
This requires $p = \exp({o(n/d_o^2)})$. If $d_o$ is bounded,
then $p$ can be as large as $e^{o({n})}$.

\citet{BraFanJia11} considered estimation as well as variable selection
and oracle properties for general concave penalties, including the lasso.
Their broader scope seems to have led to more elaborate statements and
some key conditions that are more difficult to verify than those of
Theorems~\ref{th-1} and~\ref{ThmA},
for example, their Condition 2(i) on a uniformly small spectrum
bound between $S^{(2)}(t,{\bbeta}_1)$ and its population version for a sparse
${\bbeta}_1$ in a neighborhood of ${\bbeta}^o$.

\begin{pf*}{Proof of Theorem~\ref{ThmA}}
Let\vspace*{1pt} $C_0=1$ and
$x=\lambda(\xi-1)/\{K(\xi+1)\} = \sqrt{(2/n)\log(2p/\eps)}$ in
Lemma~\ref{LemC}.
The probability of the event $|\dell({\bbeta}^o)|_{\infty}> (\xi
-1)/(\xi+1)\lambda$ is
at most $\eps$. 
The desired result follows directly from Theorem~\ref{th-1}.
\end{pf*}

\section{Compatibility and invertibility factors and restricted
eigenvalues}\label{sec4}

In Section~\ref{sec3}, the oracle inequalities in Theorems~\ref{th-1}
and~\ref{ThmA} are expressed in terms of the compatibility and weak
cone invertibility factors. However, as mentioned in the
\hyperref[sec1]{Introduction}, these quantities are still random
variables. This section provides sufficient conditions under which they
can be treated as constants. Since these factors appear in the
denominator of error bounds, it suffices to bound them from below. We
also derive a lower bound for the restricted eigenvalue to facilitate
further analysis of the Cox
model in high-dimension.
We will prove that these quantities are bounded from below by the
population version
of their certain truncated versions.

Compared with linear regression, our problem poses two additional
difficulties in the Cox model:
(a) time dependence of covariates, and (b) stochastic integration of
the Hessian over random risk sets.
Fortunately, the compatibility and weak cone invertibility factors in
Theorems~\ref{th-1} and~\ref{ThmA}
involve only the Hessian of the log-partial likelihood at the true
$\bbeta^o$, so that a martingale
argument can be used.

To simplify the statement of our results, we use
$\phi(\xi,\cO;\bSigmabar)$ to denote any of the following quantities:
%
%
\begin{equation}
\label{phi}\quad
\phi(\xi,\cO;\bSigmabar)=\kappa^2(\xi,\cO;\bSigmabar),\qquad
F_q(\xi,\cO;\bSigmabar)
\quad\mbox{or}\quad\RE^2(
\xi,\cO;\bSigmabar),
\end{equation}
where $\kappa(\xi,\cO;\bSigmabar)$, $F_q(\xi,\cO;\bSigmabar)$, and
$\RE(\xi,\cO;\bSigmabar)$ are as in (\ref{Re1}), (\ref{Fq}) and
(\ref{Re2}), respectively. If we make a claim about
$\phi(\xi,\cO;\bSigmabar)$, we mean that the claim holds for any
quantity it represents. Let $\phi_{\min}$ denote the smallest
eigenvalue. The following lemma provides some key properties of
$\phi(\xi,\cO;\bSigmabar)$ used in the derivation of its lower bounds.

%
\begin{lemma}\label{lm-RE} Let $\kappa(\xi,\cO;\bSigmabar)$,
$F_q(\xi,\cO;\bSigmabar)$,
$\RE(\xi,\cO;\bSigmabar)$ and $\phi(\xi,\cO;\bSigmabar)$ be as
in (\ref{phi}).
Let $\Sigmabar_{jk}$ be the elements of 
$\bSigmabar$ and $\bSigma$ be another nonnegative-definite matrix
with elements $\Sigma_{jk}$.

\begin{longlist}[(iii)]
\item[(i)]
For $1\le q\le2$,
\[
\qquad\min\bigl\{\kappa^2(\xi,\cO;\bSigmabar),(1+\xi)^{2/q-1}F_q(
\xi,\cO;\bSigmabar)\bigr\} \ge\RE^2(\xi,\cO;\bSigmabar) \ge
\phi_{\min}(\bSigmabar).
\]
\item[(ii)] $\phi(\xi,\cO;\bSigmabar)\ge\phi(\xi,\cO;\bSigma)
- d^o(\xi+1)^2\max_{1\le j\le k\le p}|\Sigmabar_{jk}- \Sigma
_{jk}|$.\vspace*{1pt}

\item[(iii)] If $\bSigmabar\ge\bSigma$, then $\phi(\xi,\cO;\bSigmabar
)\ge\phi(\xi,\cO;\bSigma)$.
\end{longlist}
\end{lemma}

\begin{pf}
By the H\"older inequality, $|\bb
|_q\le|\bb|_1^{2/q-1}|\bb|_2^{2-2/q}$
for all \mbox{$1 \leq q \leq 2$}.
Since $|\bb|_1\le(1+\xi
)|\bb_{\cO}|_1$
in the cone and $|\bb_{\cO}|_1\le d_o^{1/2}|\bb|_2$, we have
\[
|\bb_\cO|_1|\bb|_q/d_o^{1/q}
\le(1+\xi)^{2/q-1}|\bb_{\cO
}|_1^{2/q}|
\bb|_2^{2-2/q}/d_o^{1/q}\le(1+
\xi)^{2/q-1}|\bb|_2^2.
\]
This and $|\bb_{\cO}|_1\le d_o^{1/2}|\bb|_2$ yields part (i) by the
definition of the quantities involved.
Part (ii) follows from
$|\bb'\bSigmabar\bb- \bb'\bSigma\bb|\le|\bb|_1^2\max_{j,k}|\Sigmabar
_{jk}- \Sigma_{jk}|$
and $|\bb|_1\le(\xi+1)|\bb_{\cO}|_1\le(\xi+1)d_o^{1/q}|\bb|_q$.
Part (iii) follows immediately from the definition of the quantities
in (\ref{phi}).
\end{pf}

It follows from Lemma~\ref{lm-RE}(ii) and (iii) that quantities of
type $\phi(\xi,\cO;\bSigmabar)$ in (\ref{phi})
can be bounded from below in two ways. The first is to bound the matrix
$\bSigmabar$ from below
and the second is to approximate $\bSigmabar$ under the supreme norm
for its elements.
In the $p\gg n$ setting, our problem is essentially the rank deficiency
of $\bSigmabar$ to begin with,
so that its lower bound is still rank deficient. However, a lower bound
of the random matrix
$\bSigmabar= \ddell(\bbeta^o)$, for example, a certain truncated
version of it,
may have a smaller variability to allow an approximation by its\vadjust{\goodbreak}
population version.
This is our basic idea.
In fact, our analysis takes advantage of this argument twice to remove
different sources of randomness.

According to our plan described in the previous paragraph,
we first choose a suitable truncation of $\bSigmabar= \ddell(\bbeta
^o)$ as a lower bound of the matrix.
This is done by truncating the maximum event time under consideration.
It follows from (\ref{hessian1}) that for $t^*>0$,
%
%
\begin{equation}
\label{truncation}\quad \ddell\bigl(\bbeta^o\bigr)\ge\ddell\bigl(
\bbeta^o;t^*\bigr)\qquad\mbox{where } \ddell\bigl(\bbeta^o;t^*
\bigr) = n^{-1}\int_0^{t^*}
V_n\bigl(s,\bbeta^o\bigr)\,d{\bar N}(s).
\end{equation}
This allows us to remove the randomness from the counting process by
replacing the average counting measure $n^{-1}\,d{\bar N}(t)$ by its
compensator\break
$R_n(s,\bbeta^o)\,d\Lambda_0(s)$,
where $\Lambda_0$ is the baseline cumulative hazard function.
This approximation of $\ddell(\bbeta^o;t^*)$ can be written as
%
%
\begin{equation}
\label{hessian-2} \bSigmabar\bigl(t^*\bigr) = \int_0^{t^*}
V_n\bigl(s,\bbeta^o\bigr)R_n\bigl(s,\bbeta
^o\bigr)\,d\Lambda_0(s).
\end{equation}

To completely remove the randomness with $\bSigmabar(t^*)$,
we apply the method again by truncating the weights $e^{\bZ
_i'(t)\bbeta^o}$ with $R_n(s,\bbeta^o)$.
For $M>0$, define
%
%
\begin{equation}
\label{hessian-3} \bSigmabar\bigl(t^*;M\bigr) = \int_0^{t^*}
\hbG_n(s;M) \,d\Lambda_0(s),
\end{equation}
where $\hbG_n(t;M) \!=\! n^{-1}\sum_{i=1}^n \{\mathbf{Z}_i\!-\!
\bar{\mathbf{Z}}_n(t;M) \}^{\otimes2}
Y_i(t)\min\{M,\exp(\bZ'_i(t){\bbeta}^o)\}$ with
%
\[
\bar{\mathbf{Z}}_n(t;M)=\frac{\sum_{i=1}^n \mathbf
{Z}_i(t)Y_i(t)\min\{M,\exp(\bZ'_i(t){\bbeta}^o)\}} {
\sum_{i=1}^n Y_i(t)\min\{M,\exp(\bZ'_i(t){\bbeta}^o)\}}. 
\]
We will prove that the matrix (\ref{hessian-3}) is a lower bound of
(\ref{hessian-2}).
Suppose $\{Y_i(t),\bZ_i(t),\allowbreak t\ge0\}$ are i.i.d. stochastic processes
from $\{Y(t),\bZ(t), t\ge0\}$.
The population version of (\ref{hessian-3}) is then
%
%
\begin{equation}
\label{hessian-4} \bSigma\bigl(t^*;M\bigr) = E\int_0^{t^*}
\bG_n(s;M) \,d\Lambda_0(s),
\end{equation}
where $\bG_n(t;M) = n^{-1}\sum_{i=1}^n \{\mathbf{Z}_i-\bmu
(t;M) \}^{\otimes2}
Y_i(t)\min\{M,\exp(\bZ'_i(t){\bbeta}^o)\}$ with
\[
\bmu(t;M)=\frac{E[\mathbf{Z}(t)Y(t)\min\{M,\exp(\bZ'{\bbeta
}^o)\}]} {
E[Y(t)\min\{M,\exp(\bZ'{\bbeta}^o)\}]}. 
\]

The analysis outlined above leads to the following main result of this section.
For $\xi\ge1$ and $\cO\subset\{1,\ldots,p\}$ with $|\cO|=d_o$,
let $\phi(\xi,\cO;\bSigmabar)$ represent all quantities of interest
given in (\ref{phi}),
$\kappa(\xi,\cO)$ and $F_q(\xi,\cO)$ be the compatibility and weak
cone invertibility factors in (\ref{Fq-Cox})
with the Hessian $\ddell(\bbeta^o)$ in (\ref{hessian1}) at the true
$\bbeta$,
and $\RE(\xi,\cO;\bSigmabar)$ be as in (\ref{Re2}).
Let $L_n(t) = \sqrt{(2/n)\log t}$.

%
\begin{theorem}\label{th-RE}
$\!\!\!$Suppose $\{Y_i(t),\bZ_i(t), t\ge0\}$ are i.i.d. processes from $\{
Y(t),\allowbreak\bZ(t),t\ge0\}$ with\vadjust{\goodbreak}
$\sup_t P\{|\mathbf{Z}_i(t)-\mathbf{Z}(t)|_\infty\le K\}=P\{
\max_i N_i(\infty)\le1\}=1$.
Let $\{t^*,M\}$ be positive constants and $r_* = EY(t^*)\min\{M,\exp
(\bZ'(t^*){\bbeta}^o)\}$.
Then,
%
%
\begin{eqnarray}
\label{th-RE-1}
&& \phi\bigl(\xi,\cO;\ddell\bigl(\bbeta^o\bigr)\bigr)
\nonumber\\
&&\qquad\ge\phi\bigl(\xi,\cO;\bSigma\bigl(t^*;M\bigr)\bigr)\\
&&\qquad\quad{} - d_o(
\xi+1)^2K^2 \bigl\{ C_1L_n
\bigl(p(p+1)/\eps\bigr) + C_2 t_{n,p,\eps}^2 \bigr\}
\nonumber
\end{eqnarray}
with at least probability $1-3\eps$, where
$C_1=1+\Lambda_0(t^*)$, $C_2=(2/r_*)\Lambda_0(t^*)$
and $t_{n,p,\eps}$ is the solution of
$p(p+1) \exp\{- nt_{n,p,\eps}^2/(2+2t_{n,p,\eps}/3) \} =
\eps/2.221$.
Consequently, for $1\le q\le2$,
%
%
\begin{eqnarray}
\label{th-RE-2}
&& \min\bigl\{\kappa^2(\xi,\cO),(1+
\xi)^{2/q-1}F_q(\xi,\cO)\bigr\}
\nonumber
\\
&&\qquad\ge\RE^2\bigl(\xi,\cO;\ddell\bigl(\bbeta^o\bigr)
\bigr)
\\
&&\qquad\ge\rho_* - d_o(\xi+1)^2K^2 \bigl
\{C_1L_n\bigl(p(p+1)/\eps\bigr) + C_2
t_{n,p,\eps}^2 \bigr\}
\nonumber
\end{eqnarray}
with at least probability $1-3\eps$,
where $\rho_*=\phi_{\min}(\bSigma(t^*;M))$ with the matrix in~(\ref
{hessian-4}).
\end{theorem}

Theorem~\ref{th-RE} implies that the compatibility and cone
invertibility factors and the restricted eigenvalue
can be all treated as constants in high-dimensional Cox
model with time-dependent covariates.
We note that $C_2 t_{n,p,\eps}^2$ is of smaller order than
$L_n(p(p+1)/\eps)$ so that
the lower bounds in (\ref{th-RE-1}) and (\ref{th-RE-2}) depend on the
choice of $t^*$ and $M$
marginally through $C_1$ and $\rho_*$.
If $d^o\sqrt{(\log p)/n}$ is sufficiently small as assumed in
Theorem~\ref{ThmA}, the right-hand side of (\ref{th-RE-2}) can be
treated as $\rho_*/2$.
It is reasonable to treat $\rho_*$ as a constant since it is the
smallest eigenvalue of a
population integrated covariance matrix in (\ref{hessian-4}).

In the proof of Theorem~\ref{th-RE}, the martingale exponential
inequality in Lemma~\ref{LemC}
is used to bound the difference between (\ref{truncation}) and (\ref
{hessian-2}).
The following Bernstein inequality for $V$-statistics is used to bound
the difference
between (\ref{hessian-3}) and (\ref{hessian-4}).
This 
inequality can be viewed as an extension of the \citet{Hoe63}
inequality for sums of bounded independent variables and nondegenerate
$U$-statistics.

%
\begin{lemma}\label{lm-V-stat}
Let $X_i$ be a sequence of independent stochastic processes
and $f_{i,j}$ be functions of $X_i$ and $X_j$ with $|f_{i,j}|\le1$.
Suppose $f_{i,j}$ are degenerate in the sense of
$E[f_{i,j}|X_i] = E[f_{i,j}|X_j]=0$ for all $i\neq j$.
Let $V_n = \sum_{i=1}^n\sum_{j=1}^n f_{i,j}$. Then,
\[
P \bigl\{ \pm V_n > (nt)^2 \bigr\} \le\frac{2 \eps_n(t)(1+\eps
_n(t))}{(1+\eps_n^2(t))^2}
\le2.221 \exp\biggl(-\frac{nt^2/2}{1+t/3} \biggr),
\]
where $\eps_n(t) = e^{-(nt^2/2)/(1+t/3)}$.
\end{lemma}

Our discussion focuses on the quantities in (\ref{phi})
for the Hessian matrix $\bSigmabar= \ddell(\bbeta^o)$ evaluated at\vadjust{\goodbreak}
the true vector of coefficients.
Still, through Lemma~\ref{LemB}, Theorem~\ref{th-RE} also provides
lower bounds
for these quantities at any $\bbeta$ not far from the true $\bbeta^o$
in terms of the $\ell_1$ distance.
We formally state this result in the following corollary.

%
\begin{corollary}\label{cor-RE}
Let $\phi(\xi,\cO;\bSigmabar)$ represent any quantities in (\ref
{phi}). Then,
\[
e^{-2\eta_{\bb}}\phi\bigl(\xi,\cO;\ddell\bigl(\bbeta^o\bigr)\bigr)
\le\phi\bigl(\xi,\cO;\ddell\bigl(\bbeta^o+\bb\bigr)\bigr) \le
e^{2\eta_{\mathbf
{b}}}\phi\bigl(\xi,\cO;\ddell\bigl(\bbeta^o\bigr)\bigr),
\]
where
\[
\eta_{\mathbf{b}}=\sup_{s} 
\max_{i,j}\bigl|\mathbf{b}'
\mathbf{Z}_i(s)-\mathbf{b}'\mathbf{Z}_j(s)\bigr|.
\]
Consequently, when $|\mathbf{Z}_i(s)-\mathbf{Z}_j(s)|_\infty
\le K$,
\begin{eqnarray*}
&& \inf\bigl\{\phi\bigl(\xi,\cO;\ddell(\bbeta)\bigr)\dvtx \bigl|\bbeta
-\bbeta
^o\bigr|_1\le\eta/(2K) \bigr\}
\\
&&\qquad\ge e^{-\eta}\phi\bigl(\xi,\cO;\ddell\bigl(\bbeta^o\bigr)
\bigr)
\\
&&\qquad\ge e^{-\eta} \bigl[\rho_* - d_o(\xi+1)^2K^2
\bigl\{ C_1L_n\bigl(p(p+1)/\eps\bigr) + C_2
t_{n,p,\eps}^2 \bigr\} \bigr]
\end{eqnarray*}
under the conditions of Theorem~\ref{th-RE}.
\end{corollary}

It is worthwhile to point out that unlike typical ``small ball''
analysis based on Taylor expansion,
Corollary~\ref{cor-RE} provides nonasymptotic control of the
quantities in an $\ell_1$ ball of
constant size.
Since $\bb'\bSigmabar\bb$ appears in the numerator of the quantities
represented by $\phi(\xi,\cO;\bSigmabar)$, Corollary~\ref{cor-RE} follows
immediately from Theorem~\ref{th-RE} and (\ref{LemB2}).
It implies that the Hessian has sufficient invertibility properties in
the analysis
of the lasso when the estimator is not far from the true $\bbeta^o$ in
$\ell_1$ distance.
On the other hand, if the Hessian has sufficient invertibility
properties in a ball of fixed size,
nonasymptotic error bounds for the lasso estimator can be established.
This ``chicken and egg'' problem is directly solved in the proof of
Theorem~\ref{th-1}.



\section{Concluding remarks}\label{sec5}
This paper deals with the Cox proportional hazards regression model
when the number of time-dependent covariates $p$ is potentially much
larger than the sample size $n$.
The $\ell_1$ penalty is used to regularize the log-partial likelihood function.
Error bounds parallel to those of the lasso in linear regression are
established.
In establishing these bounds, we extend the notion of the restricted eigenvalue
and compatibility and cone invertibility factors to the Cox model.
We show that these quantities indeed provide useful error bounds.

An important issue is the choice of the penalty level $\lam$.
Theorem~\ref{ThmA} requires a $\lam$ slightly larger than $K\sqrt
{(2/n)\log p}$,
where $K$ is a uniform upper bound for the range of individual real covariates.
This indicates that the lasso is tuning insensitive since the
theoretical choice does not
depend on the unknowns.
In practice, cross validation can be used to fine tune the penalty
level $\lam$.
Theoretical investigation of the performance of the lasso with
cross-validated $\lam$,
an interesting and challenging problem in and of itself even in the
simpler linear regression model,
is beyond the scope of this paper.\vadjust{\goodbreak}\looseness=-1

General concave penalized estimators in the Cox regression model have been
considered in \citet{BraFanJia11}
where oracle inequalities and properties of certain local solutions are
considered.
\citet{ZhaZha12} has provided a unified treatment of global and
local solutions
for concave penalized least squares estimators in linear regression.
Since this unified
treatment relies on an oracle inequality for the global solution based
on the cone
invertibility factor, the results in this paper point to a possible
extension of such a unified
treatment of global and local solutions of general concave regularized methods
in the Cox regression model.

%
\begin{appendix}\label{app}
\section*{Appendix}
Here we prove Lemmas~\ref{LemH1},~\ref{LemB},~\ref{LemC} and~\ref{lm-V-stat}
and Theorems~\ref{th-1} and~\ref{th-RE}.

\begin{pf*}{Proof of Lemma~\ref{LemH1}} 
Since $\ell({\bbeta})$ is a convex function,
$D^s(\bhbeta,\bbeta) =\tbtheta'\{\dell
({\bbeta}^o+\tbtheta) - \dell({\bbeta}^o)\}\ge0$,
so that the first inequality holds.
Since $\ttheta_j=\hbeta_j$ for $j\in\cO^c$, (\ref{KKT}) gives
\begin{eqnarray*}
\tilde{\btheta} \bigl\{\dell\bigl(\bbeta^o+\tilde{\btheta} \bigr)
- \dell\bigl(\bbeta^o\bigr) \bigr\}
&=& \sum_{j\in\cO^c}{\ttheta}_j \bigl(\dell
\bigl(\bbeta^o+\tilde{\btheta}\bigr) \bigr)_j + \sum
_{j\in\cO}{\ttheta}_j \bigl(\dell\bigl(
\bbeta^o+\tilde{\btheta}\bigr) \bigr)_j+\tilde{
\btheta}' \bigl(-\dell\bigl(\bbeta^o\bigr) \bigr)\hspace*{5pt}
\\[-2pt]
&\le&\sum_{j\in\cO^c}\hbeta_j \bigl(-\lam
\sgn(\hbeta_j) \bigr) + \sum_{j\in\cO}|{
\ttheta}_j|\lam+ |\tilde{\btheta}|_1 z^*\hspace*{5pt}
\\[-2pt]
&=& \sum_{j\in\cO^c} - \lam|\tilde{\theta}_j|
+ |\tilde{\btheta}_{\cO}|_1\lambda+ z^* |\tilde{
\btheta}_{\cO}|_1 + z^* |\tilde{\btheta}_{\cO^c}|_1\hspace*{5pt}
\\[-2pt]
&=& \bigl(z^*-\lambda\bigr)|\tilde{\btheta}_{\cO^c}|_1 +
\bigl(\lambda+z^*\bigr)|\tilde{\btheta}_{\cO}|_1.\hspace*{5pt}
\end{eqnarray*}
%
Thus the second inequality in (\ref{LemH1E1}) holds.
Note that the inequality in the third line above requires
$ (\dell(\bbeta^o+\tilde{\btheta}) )_j = -\lam\sgn(\hbeta_j)$
only in the set $\cO^c\cap\{j\dvtx \hbeta_j\neq0\}$, since
$\tilde{\theta}_j = \hat{\beta}_j - \beta^o_j = 0$ when $j \in
\mathcal{O}^c$ and $\hat{\beta}_j = 0$.
\end{pf*}

\begin{pf*}{Proof of Lemma~\ref{LemB}}
We use similar notation as in \citet{HjoPol93}. Let $a_i = a_i(s)
= \mathbf{b}'\{\mathbf{Z}_i(s)-\bar{\mathbf{Z}}_n(s,{\bbeta})\}$,
$w_{i}=w_i(s)=Y_i(s)\exp[{\bbeta}'\mathbf{Z}_i(s)]$ and $c=c(s) =
(\max_ia_i(s)+\min_ia_i(s))/2$. Clearly,
$\max_i|a_i-c|\le(1/2)\eta_{\mathbf{b}}$. By the definition of
$\bar{\mathbf{Z}}_n(t,{\bbeta})$,
\begin{eqnarray*}
&& \mathbf{b}'\bigl\{\bar{\mathbf{Z}}_n(s,{\bbeta}+
\mathbf{b}) - \bar{\mathbf{Z}}_n(s,{\bbeta})\bigr\}
\\[-2pt]
&&\qquad= \sum_i \mathbf{b}'
\mathbf{Z}_i(s) w_i e^{\mathbf
{b}'\mathbf{Z}_i(s)} \Big/\sum
_i w_i e^{\mathbf{b}'\mathbf{Z}_i(s)} - \sum
_i \mathbf{b}'\mathbf{Z}_i(s)
w_i \Big/\sum_i w_i
\\[-2pt]
&&\qquad= \sum_i a_i w_i
e^{a_i} \Big/\sum_i w_i
e^{a_i} - \sum_i a_i
w_i \Big/\sum_i w_i
\\[-2pt]
&&\qquad= \sum_{i,j} w_{i} w_{j}
a_i\bigl(e^{ a_i} - e^{ a_j}\bigr) \Big/\sum
_{i,j} w_{i} w_{j} e^{ a_i}
\\[-2pt]
&&\qquad= \sum_{i,j} w_{i} w_{j}
(a_i-a_j) \bigl(e^{ a_i-c} - e^{ a_j-c}
\bigr) \Big/\sum_{i,j} 2w_{i}
w_{j} e^{ a_i-c}
\\[-2pt]
&&\qquad\ge\exp\Bigl({-2\max_i}|a_i-c| \Bigr)\sum
_{i,j} w_{i} w_{j}
(a_i-a_j)^2 \Big/\sum
_{i,j} 2w_{i} w_{j}
\\[-2pt]
&&\qquad\ge\exp(-\eta_{\mathbf{b}} )\sum_{i}
w_{i} a_i^2 \Big/\sum
_{i} w_{i},
\end{eqnarray*}
where the first inequality comes from $(e^y-e^x)/(y-x)\ge e^{-(|y|\vee
|x|)}$ and, since $\sum_iw_ia_i=0$,
the second one from $\sum_{i,j} w_{i} w_{j} (a_i-a_j)^2=2\sum_{i}
w_{i} \sum_{i} w_{i} a_i^2$.
Thus, since $a_i^2=\mathbf{b}'\{\mathbf{Z}_i(s)-\bar{\mathbf
{Z}}_n(s,{\bbeta})\}^{\otimes2}\mathbf{b}$,
(\ref{hessian1}) and (\ref{gradient1}) give
\[
e^{-\eta_{\mathbf{b}}}\mathbf{b}'\ddell({\bbeta})\mathbf{b} =
\frac{e^{-\eta_{\mathbf{b}}}}{n} \int_0^\infty\sum
_{i=1}^n w_{i} a_i^2
\Biggl(\sum_{i=1}^n w_{i}
\Biggr)^{-1} \,d\Nbar(s) \le\mathbf{b}'\bigl\{\dell({
\bbeta}+\mathbf{b})-\dell({\bbeta})\bigr\}.
\]
This implies the lower bound in (\ref{LemB1}). Similarly, the lower
bound in (\ref{LemB2}) follows from
%
\begin{eqnarray*}
\ddell(\bbeta+ \bb) &=& \frac{1}{n} \int_0^\infty
\frac{\sum_{i,j} w_iw_j\{\bZ_i(s)\bZ
_i'(s) - \bZ_i(s)\bZ_j'(s)\}e^{a_i+a_j}} {
\sum_{i,j} w_iw_je^{a_i+a_j}}\,d\Nbar(s)
\\
&=& \frac{1}{n} \int_0^\infty
\frac{\sum_{i,j} w_iw_j (\bZ
_i(s)-\bZ_j(s) )^{\otimes2}e^{(a_i-c)+(a_j-c)}}{\sum_{i,j}
2w_iw_je^{(a_i-c)+(a_j-c)}}\,d\Nbar(s)
\end{eqnarray*}
and
\[
\ddell(\bbeta) = \frac{1}{n} \int_0^\infty
\frac{\sum_{i,j}w_i w_j (\bZ
_i(s)-\bZ_j(s) )^{\otimes2}}{\sum_{i,j} 2w_i w_j}\,d\Nbar(s).
\]
The proof of the upper bounds in (\ref{LemB1}) and (\ref{LemB2}),
nearly identical to the proof of the lower bounds, is omitted.
\end{pf*}

\begin{pf*}{Proof of Lemma~\ref{LemC}}
Applying the union bound and changing the scale of the covariates if necessary,
we assume without loss of generality that $p=K=1$.
In this case
%
\[
\dell\bigl({\bbeta}^o\bigr) = \frac{1}{n}\sum
_{i=1}^n\int_0^\infty
a_i(s)\,dN_i(s) = \frac{1}{n}\sum
_{i=1}^n\int_0^\infty
a_i(s)\,dM_i(s),
\]
where $a_i(t)={Z}_{i1}(t)- \bar{Z}_{n,1}(t)$, $i=1,\ldots, n$, are
predictable and
satisfy \mbox{$|a_i(t)|\le1$}. Thus, (\ref{lm-C-1}) follows from (\ref{lm-C-0}).

Let $t_j$ be the time of the $j$th jump of the process
$\sum_{i=1}^n \int_0^\infty Y_i(t) \,dN_i(t)$, $j=1,\ldots,m$ and $t_0=0$.
Then, $t_j$ are stopping times.
For $j=0,\ldots,m$, define
%
%
\begin{equation}
\label{martX} X_j=\sum_{i=1}^n
\int_0^{t_j} a_i(s)\,dN_i(s)
=\sum_{i=1}^n\int_0^{t_j}
a_i(s)\,dM_i(s).
\end{equation}
Since $M_i(s)$ are martingales and $a_i(s)$ are predictable,
$\{X_j, j=0,1,\ldots\}$ is a martingale with the difference
$|X_j-X_{j-1}|\le\max_{s,i}|a_i(s)|\le1$.
Let $m$ be the greatest integer lower bound of $C_0^2n$.
By the martingale version of the \citet{Hoe63} inequality
[\citet{Azu67}],
\[
\rP\bigl(|X_m|>nC_0x\bigr)\le2\exp\bigl(-n^2C_0^2x^2/(2m)
\bigr) \le 2e^{-nx^2/2}.
\]
By (\ref{martX}), $X_m=n\dell({\bbeta}^o)$ if and only if
$\sum_{i=1}^n \int_0^\infty Y_i(t) \,dN_i(t)\le m$.
Thus, the left-hand side of (\ref{lm-C-1}) is no greater than
$\rP(|X_m|>nC_0x)\le 2e^{-nx^2/2}$.
\end{pf*}

\begin{pf*}{Proof of Lemma~\ref{lm-V-stat}}
For integers $j,m,i_1,\ldots,i_m$, let $\#(j;i_1,\ldots,i_m)$ be the
number of appearances of $j$ in the sequence $\{i_1,\ldots,i_m\}$.
Since $f_{i,j}$ are degenerate,
\begin{eqnarray*}
E(\pm V_n)^m &=& \sum_{1\le i_1,\ldots,i_{2m}\le n}
(\pm1)^m f_{i_1,i_2}\cdots f_{i_{2m-1},i_{2m}}
\\
&\le& \sum_{1\le i_1,\ldots,i_{2m}\le n} \prod
_{j=1}^n I\bigl\{\# (j;i_1,\ldots,i_{2m})\neq 1\bigr\}.
\end{eqnarray*}
This is due to the fact that all terms with exactly one appearance of
an index $j$ have
zero expectation and all other terms are bounded by 1.
Let $E_0$ be the expectation under which $i_1,\ldots,i_{2m}$ are
i.i.d. uniform variables in $\{1,\ldots,n\}$
and $k_j = \#(j;i_1,\ldots,i_{2m})$.
Since $(k_1,\ldots,k_n)$ is multinomial$(2m,1/n,\ldots,1/n)$,
the above inequality can be written as
\[
E(\pm V_n)^m \le n^{2m} E_0\prod
_{j=1}^n I\{k_j\neq 1\} = (2m)!
\sum_{k_1+\cdots+k_n=2m} \prod_{j=1}^n
\frac{I\{k_j\neq 1\}}{k_j!}.
\]
Let $f_0(x) = \sum_{m=0}^\infty x^m/(2m)! = \cosh(|x|^{1/2})I\{x\ge
0\} + \cos(|x|^{1/2})I\{x<0\}$
and $\lam= t/(1+t/3)$. It follows from the above moment inequality that
\begin{eqnarray*}
Ef_0\bigl(\pm\lam^2 V_n\bigr) &=& \sum
_{m=0}^\infty\lam^{2m} E(\pm
V_n)^m/(2m)!
\\
&\le& \sum_{m=0}^\infty\lam^{2m}
\sum_{k_1+\cdots+k_n=2m} \prod_{j=1}^n
\frac{I\{k_j\neq 1\}}{k_j!}
\\
&\le& \sum_{m=0}^\infty\lam^{m}
\sum_{k_1+\cdots+k_n=m} \prod_{j=1}^n
\frac{I\{k_j\neq 1\}}{k_j!}
\\
&=& \Biggl(\sum_{k=0}^\infty
\lam^k I\{k\neq1\}/k! \Biggr)^n.
\end{eqnarray*}
Since\vspace*{1pt} $\sum_{k=0}^\infty\lam^k I\{k\neq1\}/k! \le1+(\lam
^2/2)/(1-\lam/3)=1+\lam t/2$, we find\break
$Ef_0(\pm\lam^2 V_n) \le e^{n\lam t/2}$.\vadjust{\goodbreak}
Consequently, the monotonicity of $f(x) = \cosh(x^{1/2})$ for $x>0$
and the lower bound
$f(x)\ge-1$ allow us to apply the Markov inequality as follows:
\begin{eqnarray*}
P \bigl\{\pm V_n > (nt)^2 \bigr\} &\le& P \bigl\{1+
f_0\bigl(\pm\lam^2 V_n\bigr) > 1+
f_0\bigl((\lam nt)^2\bigr) \bigr\}
\\
&\le& \bigl\{1+f_0\bigl((n\lam t)^2\bigr)\bigr
\}^{-1}E\bigl\{1+f_0\bigl(\pm\lam^2
V_n\bigr)\bigr\}
\\
&\le& \bigl\{1+\cosh(n\lam t)\bigr\}^{-1}\bigl(1+e^{n\lam t/2}
\bigr) 
\\
& = & 2 e^{- n\lam t/2}\bigl(1+e^{-n\lam t/2}\bigr)/\bigl(1+e^{-n\lam t}
\bigr)^2.
\end{eqnarray*}
The conclusion follows from $2\max_{0\le x\le1}(1+x)/(1+x^2)^2\le
2.221$.
\end{pf*}

\begin{pf*}{Proof of Theorem~\ref{th-1}}
Let $\tbtheta=\bhbeta-{\bbeta}^o\neq0$
and $\mathbf{b}=\tbtheta/|\tbtheta|_1$.
It follows from the convexity of $\ell({\bbeta}^o+x\mathbf{b})$,
as a function of $x$, and
Lemma~\ref{LemH1} that, in the event
$|\dell({\bbeta}^o)|_{\infty}\le(\xi-1)/(\xi+1)\lam$,
%
%
\begin{equation}
\label{pf-th-1-1} \mathbf{b}'\bigl\{\dell\bigl({
\bbeta}^o+x\mathbf{b}\bigr) - \dell\bigl({\bbeta}^o\bigr)
\bigr\} + \frac{2\lam}{\xi+1} |\mathbf{b}_{\cO^c}|_1 \le
\frac{2\xi\lam}{\xi+1} |\mathbf{b}_{\cO}|_1
\end{equation}
for $ x\in[0, |\tbtheta|_1]$ and $\mathbf{b} \in
\mathscr{C}(\xi,\cO)$.
Consider all nonnegative $x$ satisfying (\ref{pf-th-1-1}). We need to
establish a lower bound for
\[
\mathbf{b}'\bigl\{\dell\bigl({\bbeta}^o+x\mathbf{b}
\bigr)-\dell\bigl({\bbeta}^o\bigr)\bigr\} = \frac{1}{n}\int
_0^\infty\mathbf{b}'\bigl\{
\bar{Z}_n\bigl(s,{\bbeta}^o+x \mathbf{b}\bigr) -
\bar{Z}_n\bigl(s,{\bbeta}^o\bigr)\bigr\}\,d\Nbar(s).
\]
Since $\eta_{x\mathbf{b}} = \max_{0\le s\le1}\max_{i,j}|x\mathbf{b}'
\mathbf{Z}_i(s)-x\mathbf{b}'\mathbf{Z}_j(s)| \le
Kx|\mathbf{b}|_1=Kx$,
Lem\-ma~\ref{LemB} yields
%
%
\begin{eqnarray}
\label{pf-th-1-2a}
x\mathbf{b}'\bigl\{\dell\bigl({
\bbeta}^o+x\mathbf{b}\bigr)-\dell\bigl({\bbeta}^o\bigr)
\bigr\} &\ge& x^2\exp(-\eta_{x\mathbf{b}} ) \mathbf
{b}'\ddell\bigl({\bbeta}^o\bigr)\mathbf{b} \nonumber\\[-8pt]\\[-8pt]
&\ge&
x^2\exp(-Kx) \mathbf{b}'\ddell\bigl({
\bbeta}^o\bigr)\mathbf{b}.\nonumber
\end{eqnarray}
This, combined with (\ref{pf-th-1-1}) and the definition of $\kappa
(\xi,\cO)$, gives
\begin{eqnarray*}
x e^{-Kx}\kappa^2(\xi,\cO)|\mathbf{b}_{\cO}|_1^2
/d_o &\le& x e^{-Kx} \mathbf{b}'\ddell\bigl({
\bbeta}^o\bigr)\mathbf{b}
\\
&\le& \frac{2\xi\lam}{\xi+1} |\mathbf{b}_{\cO}|_1 -
\frac
{2\lam}{\xi+1}|\mathbf{b}_{\cO^c}|_1
\\
&=& 2\lam|\mathbf{b}_{\cO}|_1 - \frac{2\lam}{\xi+1}
\\
&\le&\lam(\xi+1)|\mathbf{b}_{\cO}|_1^2/2.
\end{eqnarray*}
In other words, any $x$ satisfying (\ref{pf-th-1-1}) must satisfy
%
%
\begin{equation}
\label{pf-th-1-2} Kx\exp(-Kx)\le\frac{K(\xi+1)\lam
d_o}{2\kappa^2(\xi,\cO)} = \tau.
\end{equation}
Since $\bb'\{\dell({\bbeta}^o+x\mathbf{b}) - \dell({\bbeta
}^o)\}$ is an increasing
function of $x$ due to the convexity of~$\ell$, the set of all
nonnegative $x$
satisfying (\ref{pf-th-1-1}) is a closed interval $[0, \tx]$ for some
$\tx>0$.
Thus,\vadjust{\goodbreak} (\ref{pf-th-1-2}) implies $K\tx\le\eta$,
the smaller solution of $\eta e^{-\eta}=\tau$.
This yields
\[
|\tbtheta|_1\le\tx\le\frac{\eta}{K} =\frac{e^{\eta}\tau}{K}=
\frac{e^{\eta}(\xi+1)\lam d_o}{2\kappa
^2(\xi,\cO)}
\]
in (\ref{th-1-1}).
The first part of (\ref{th-1-1}) follows from (\ref{Re1}), (\ref
{Fq-Cox}), (\ref{LemB1}) and (\ref{LemH1E1}), due to
%
\[
e^{-\eta} \kappa^2(\xi,\cO)|\tbtheta_{\cO}|_1^2/d^o
\le e^{-\eta} \tbtheta'\ddell\bigl({\bbeta}^o
\bigr) \tbtheta\le D^s(\bhbeta,\bbeta) \le\frac{2\xi\lam|\tbtheta_{\cO
}|_1}{\xi+1}.
\]

Finally, it follows from the definition of $F_q(\xi,\cO)$, (\ref
{pf-th-1-2a}) and
(\ref{pf-th-1-1}) that, for $x = |\tbtheta|_1$,
\begin{eqnarray*}
x e^{-\eta}&\le&\frac{x e^{-Kx} \bb'\ddell({\bbeta}^o)\mathbf{b}} {
F_q(\xi,\cO)(|\mathbf{b}_{\cO}|_1 /d_o^{1/q})|\mathbf{b}|_q} \le\frac
{\mathbf{b}'\{\dell({\bbeta}^o+x\mathbf{b})-\dell
({\bbeta}^o)\}} {
F_q(\xi,\cO)(|\mathbf{b}_{\cO}|_1 /d_o^{1/q})|\mathbf{b}|_q} \\[-2pt]
&\le&
\frac{2\xi\lam d_o^{1/q}}{(\xi+1)F_q(\xi,\cO)|\bb|_q}.
\end{eqnarray*}
This gives the second inequality in (\ref{th-1-2}) due to
$|\bhbeta-{\bbeta}^o|_q=|\tbtheta|_1|\mathbf{b}|_q$.
\end{pf*}

\begin{pf*}{Proof of Theorem~\ref{th-RE}}
Let
\begin{eqnarray*}
a(s)&=&\bigl(V_n\bigl(s,\bbeta^o\bigr)
\bigr)_{jk}/K^2 \\[-2pt]
&=&\sum_{i=1}^nw_{ni}
\bigl(t,\bbeta^o\bigr)\bigl\{Z_{i,j}(s)-{\bar
Z}_{n,j}(s)\bigr\} \bigl\{ Z_{i,k}(s)-{\bar
Z}_{n,k}(s)\bigr\}/K^2.
\end{eqnarray*}
It follows from Lemma~\ref{LemC}(i) with $a_i(s)=a(s)$ and $C_0=1$ that
\begin{eqnarray*}
&&
P \biggl\{ \biggl| \biggl(\int_0^{t^*}V_n
\bigl(s,\bbeta^o\bigr)\,d{\bar N}(s) - \int_0^{t^*}V_n
\bigl(s,\bbeta^o\bigr)R_n\bigl(s,\bbeta^o
\bigr)\,d\Lambda_0(s) \biggr)_{jk} \biggr| > K^2 x
\biggr\}\\[-2pt]
&& \le2e^{-nx^2/2}.
\end{eqnarray*}
Thus, $P\{\max_{j,k}|(\ddell(\bbeta^o;t^*) - \bSigmabar
(t^*)|_{j,k}\ge K^2L_n(p(p+1)/\eps)\}\le\eps$
by the union bound and the respective definitions of $\ddell(\bbeta
^o;t^*)$ and $\bSigmabar(t^*)$ in
(\ref{truncation}) and (\ref{hessian-2}).
Consequently, by (\ref{truncation}) and Lemma~\ref{lm-RE}(iii) and (ii)
%
%
\begin{eqnarray}
\label{pf-th-RE-1}\quad
&&
P \bigl\{\phi\bigl(\xi,\cO;\ddell\bigl(\bbeta^o
\bigr)\bigr) \ge\phi\bigl(\xi,\cO;\bSigmabar\bigl(t^*\bigr)\bigr) -
d_o(\xi+1)^2K^2L_n\bigl(p(p+1)/
\eps\bigr) \bigr\}\nonumber\\[-8pt]\\[-8pt]
&&\qquad\ge1-\eps.\nonumber
\end{eqnarray}
Let us take the sample mean of $i$-indexed quantities with
weights\break
$Y_i(t)\min\{M,\allowbreak\exp(\bZ'_i(t){\bbeta}^o)\}$, so that
$\bar{\mathbf{Z}}_n(t;M)$ is the sample mean of $\bZ_i(t)$.
Since $V_n(t,\bbeta^o)R_n(t,\allowbreak\bbeta^o) = \hbG_n(t;\infty)$,
\begin{eqnarray*}
\bu' \hbG_n(t;\infty)\bu&\ge&\frac{1}{n}\sum
_{i=1}^n \bigl[\bu' \bigl\{
\mathbf{Z}_i- \bar{\mathbf{Z}}_n(t;\infty) \bigr\}
\bigr]^2 Y_i(t)\min\bigl\{M,\exp\bigl(
\bZ'_i(t){\bbeta}^o\bigr)\bigr\} \\[-2pt]
&\ge&
\bu'\hbG_n(t;M)\bu.
\end{eqnarray*}
Thus, by the definition of $\bSigmabar(t^*;M)$ in (\ref{hessian-3})
and Lemma~\ref{lm-RE}(iii),
%
%
\begin{equation}
\label{pf-th-RE-2} \phi\bigl(\xi,\cO;\bSigmabar\bigl(t^*\bigr)\bigr) \ge
\phi
\bigl(\xi,\cO;\bSigmabar\bigl(t^*;M\bigr)\bigr).
\end{equation}
In addition, the relationship between the sample second moment and
variance gives
\[
\bG_n(t;M) = \hbG_n(t;M) + \bigl\{\bar{
\mathbf{Z}}_n(t;M) - \bmu(t;M) \bigr\}^{\otimes2}
\]
by the definition of $\bG_n(t;M)$ and $\hbG_n(t;M)$, so that (\ref
{hessian-3}) can be written as
%
%
\begin{eqnarray}
\label{pf-th-RE-3}
\bSigmabar\bigl(t^*;M\bigr) &=& \int_0^{t^*}
\bG_n(s;M)\,d\Lambda_0(s)\nonumber\\[-8pt]\\[-8pt]
&&{} - \int_0^{t^*}
\bigl\{\bar{\mathbf{Z}}_n(t;M) - \bmu(t;M) \bigr\}^{\otimes
2}\,d
\Lambda_0(s).\nonumber
\end{eqnarray}

We first bound the second term on the right-hand side of (\ref
{pf-th-RE-3}). Define
\begin{eqnarray*}
R_n(t;M)&=& \frac{1}{n}\sum_{i=1}^n
Y_i(t) \min\bigl\{M,\exp\bigl(\bZ'_i(t){
\bbeta}^o\bigr)\bigr\},
\\
\bDelta(t;M) &=& R_n(t;M)\bigl\{\bar{\mathbf{Z}}_n(t;M)
- \bmu(t;M)\bigr\}
\\
&=& \frac{1}{n}\sum_{i=1}^n
Y_i(t) \min\bigl\{M,\exp\bigl(\bZ'_i(t){
\bbeta}^o\bigr)\bigr\} \bigl\{\bZ_i(t) - \bmu(t;M)\bigr
\}.
\end{eqnarray*}
Since $Y_i(t)$ is nonincreasing in $t$,
%
%
\begin{equation}
\label{pf-th-RE-4}\quad 0\le\int_0^{t^*} \bigl\{\bar{
\mathbf{Z}}_n(t;M) - \bmu(t;M) \bigr\} ^{\otimes2}\,d
\Lambda_0(s) \le\frac{\int_0^{t^*}\bDelta^{\otimes2}(t;M)\,d\Lambda
_0(s)}{R_n^{2}(t^*,M)}.
\end{equation}
Since $R_n(t^*,M)$ is the average of i.i.d. variables uniformly bounded
by $M$ and $ER_n(t^*,M) = r_*$,
the \citet{Hoe63} inequality gives
\[
P \bigl\{R_n\bigl(t^*,M\bigr) < r_*/2 \bigr\} \le e^{- n r_*^2/(8M^2)}.
\]
Since $\bDelta(t;M)$ is an average of i.i.d. mean zero vectors,
\[
\biggl(n^2\int_0^{t^*}\bDelta^{\otimes2}(t;M)\,d\Lambda_0(s)
\biggr)_{jk}
\]
is a degenerate $V$-statistic
for each $(j,k)$. Moreover, since the summands of these $V$-statistics
are all bounded by
$K^2\Lambda_0(t^*)$, Lemma~\ref{lm-V-stat} yields
\begin{eqnarray*}
&&
\max_{1\le j,k\le p} P \biggl\{\pm\biggl(\int_0^{t^*}
\bDelta^{\otimes2}(t;M)\,d\Lambda_0(s) \biggr)_{jk} >
K^2\Lambda_0\bigl(t^*\bigr) t^2 \biggr\}\\
&&\qquad\le2.221 \exp\biggl(\frac{- nt^2/2}{1+t/3} \biggr).
\end{eqnarray*}
Thus, by (\ref{pf-th-RE-3}), (\ref{pf-th-RE-4}), the above two
probability bounds and Lemma~\ref{lm-RE}(ii),
%
%
\begin{eqnarray}
\label{pf-th-RE-5} && \phi\bigl(\xi,\cO;\bSigmabar\bigl(t^*;M\bigr)\bigr)
\nonumber\\
&&\qquad\ge\phi\biggl(\xi,\cO;\int_0^{t^*}
\bG_n(s;M)\,d\Lambda_0(s) \biggr) \\
&&\qquad\quad{}- d^o(
\xi+1)^2K^2\Lambda_0\bigl(t^*\bigr)
t_{n,p,\eps}^2/(r_*/2)
\nonumber
\end{eqnarray}
with at least probability $1-e^{- n r_*^2/(8M^2)} - \eps$.

Finally, by (\ref{hessian-4}), $\int_0^{t^*}\bG_n(s;M)\,d\Lambda
_0(s)$ is an average of i.i.d. matrices with mean
$\bSigma(t^*;M)$ and the summands of $ (\int_0^{t^*}\bG
_n(s;M)\,d\Lambda_0(s) )_{jk}$ are
uniformly bounded by $K^2\Lambda_0(t^*)$, so that the \citet{Hoe63}
inequality gives
\begin{eqnarray*}
&&
P \biggl\{\max_{j,k} \biggl| \biggl(\int_0^{t^*}
\bG_n(s;M)\,d\Lambda_0(s) - \bSigma\bigl(t^*,M\bigr)
\biggr)_{jk} \biggr| > K^2\Lambda_0\bigl(t^*\bigr) t
\biggr\} \\
&&\qquad\le p(p+1) e^{- n t^2/2}.
\end{eqnarray*}
By (\ref{pf-th-RE-1}), (\ref{pf-th-RE-2}), (\ref{pf-th-RE-5}), the
above inequality with $t=L_n(p(p+1)/\eps)$
and Lem\-ma~\ref{lm-RE}(ii),
\begin{eqnarray*}
&&
\phi\bigl(\xi,\cO;\ddell\bigl(\bbeta^o\bigr)\bigr) 
\\
&&\qquad\ge \phi\biggl(\xi,\cO;\int
_0^{t^*}\bG_n(s;M)\,d
\Lambda_0(s) \biggr)
\\
&&\qquad\quad{} - d_o(\xi+1)^2K^2 \bigl\{L_n
\bigl(p(p+1)/\eps\bigr) + (2/r_*)\Lambda_0\bigl(t^*\bigr)
t_{n,p,\eps}^2 \bigr\}
\\
&&\qquad\ge \phi\bigl(\xi,\cO;\bSigma\bigl(t^*,M\bigr) \bigr)
\\
&&\qquad\quad{} - d_o(\xi+1)^2K^2 \bigl\{\bigl(1+
\Lambda_0\bigl(t^*\bigr)\bigr)L_n\bigl(p(p+1)/\eps\bigr)
+ (2/r_*)\Lambda_0\bigl(t^*\bigr) t_{n,p,\eps}^2
\bigr\}
\end{eqnarray*}
with at least probability $1-e^{- n r_*^2/(8M^2)} - 3\eps$.
Since
\[
\phi\bigl(\xi,\cO;\bSigma\bigl(t^*,M\bigr) \bigr)\ge\RE^2 \bigl(
\xi,\cO;\bSigma\bigl(t^*,M\bigr) \bigr)\ge\rho_*
\]
by Lemma~\ref{lm-RE}(i) and the definition in (\ref{Re2}),
the conclusion follows.
\end{pf*}
\end{appendix}

\section*{Acknowledgments}

We are grateful to two anonymous reviewers, the Associate Editor and
Editor for their helpful comments which led to considerable
improvements in the paper. We also wish to thank a reviewer for
bringing to our attention the work of \citet{GafGui12},
\citet{Lem} and \citet{KonNan} during the revision process of
this paper.


%

\printaddresses

\end{document}